\documentclass[10pt]{article}
\advance\textheight by\topskip
\def\title#1{\centerline{\LARGE\bf #1}\vskip .5em}
\def\author#1{\noindent{#1}\vskip .5em}
\def\address#1{\noindent{\it #1}\vskip .5em}

\usepackage[all]{xy} \SelectTips{cm}{} 
\def\tsk{\ .}
\def\kbl{\ ,}
\let\mzl \leq
\def\inv{^{\mit -1}}
\let\trik\bigtriangleup
\let\dot\bullet

\def\sqr#1#2{{\vcenter{\vbox{\hrule height .#2pt
\hbox{\vrule width .#2 pt height #1 pt \kern #1pt
\vrule width .#2 pt}
\hrule height .#2pt}}}}
\def\square{\sqr 64}

\def\irodymopabaiga{\vskip 1pt
\hfill $\square$\vskip \baselineskip}
\begin{document}

\title{Notion of a virtual derivative}
\author{Gintaras VALIUKEVI\v CIUS}
\address{Institute of Mathematics and Informatics, Vilnius University,
 Akademijos 4, LT-08663, Vilnius, Lithuania\\
e-mail: geva@ktl.mii.lt}

{\bf Abstract.} Formal graphical procedures to calculate function's derivative 
 are proposed. This can be applied to calculate expressions of
geometric objects, construct approximate schemes in numerical analysis
of solution of differential equation.
\vskip 12pt
{\bf KEY WORDS:} derivative, graphs, virtual graphs, ordering
of virtual graphs, symmetry coefficient.
\vskip 12 pt
{\bf 1. Introduction.}
Finding derivatives of functions becomes a fatal obstacle in more complex
analysis. Partial derivatives of a function over Euclidean space were
used in Riemann's geometry from the 19-th century [1]. There was an attempt  to
escape such  cumbersome expressions with the help of notions from algebraic
geometry [2]. We can call such method as an algebraic specification of geometric
objects. It has its own difficulties. V. Arnold freely applied such abstract
notions in the analysis of infinite dimensional spaces and described the flow
of an ideal fluid governed by Euler equations as geodesic in the group of
diffeomorphisms [3]. To calculate concrete derivatives  for a detailed
study we need an easy notation of derivatives. Everyone who tryed  to define
concrete topologies in some infinite dimensional spaces will agree with this.
Specialists of algebraic geometry like S. Lang knew the shortages
of algebraic specifications and suggested  the derivative's notion 
for a function over a Banach space [4]. These new derivatives needed a new more
precise marking. That is obvious for derivatives of higher order. I started
to mark derivatives by graphs, and applied this to calculate derivatives for
the Runge Kutta algorithm [5]. Later I applied the graphs of derivatives to
describe asymptotic expansions in the Banach algebra of probabilities on
a lattice [6]. In this work we present an overview of methods to
mark derivatives, define  more explicitly the natural ordering of graphs.
Graph constructions for other new geometric or asymptotic problems remain
to be found.

{\bf 2. Virtual derivative.}
\newdimen \zingsnis\zingsnis=24pt
\newdimen\aukstis\aukstis=36pt
\newdimen\grynas\grynas=0.0pt
Every mapping can be pictured as an arrow:
$$\xy *{\bullet} ,c+<0cm,\aukstis>*{\bullet}**\dir{-}?>*\dir{>}\endxy$$
The target of the arrows will be upper directed, therefore we shall
use the arrows without heads.

A mapping of several arguments will be pictured by means of the
corresponding number of arrows having the common vertex:

\vskip 12pt
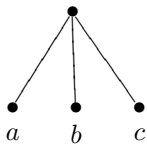
$$\xy 
@={\save *{\bullet}
\afterPOS{\POS +<\grynas, -10 pt>*{a}},\restore\POS
+<\zingsnis,\grynas>\save *{\bullet}
\afterPOS{\POS +<\grynas, -10pt>*{b}},\restore\POS
+<\zingsnis,\grynas> *{\bullet}\save
\afterPOS{\POS +<\grynas,-10pt>*{c}}} **@{}?<>(.5),
+<0cm,\aukstis>="virsus"*{\bullet}
@@{;"virsus"**@{}?(.95)**\dir{-}?>}
\endxy$$
\vskip 12pt

We choose an order of arguments pictured in the graph, 
i. e. we know which argument is the
 first, second, and so on. The composition of mappings will be pictured as
 a graph with wedges corresponding to the mappings under consideration.
 We get a tree with the top point marking the target space, and  the roots
 marking the arguments in the source spaces. Every  root  
will be called { \it an entrance of the graph}. 
The top vertice would be {\it an outlet of the graph}.
The {\it youngest  mappings} start
 with entrances, and the {\it eldest mapping} 
ends with the top vertex of the tree.
We pretend a composition of mappings to  call as a {\it joint mapping}.

If each mapping with $n$ arguments is  written as $\langle x_1, x_2, \dots,
x_n\rangle\circ f$, then the tree of joint mapping can be writen as
$$\hskip-10 pt\langle\langle x_1, \dots, x_{n1}\rangle\circ f1,
\langle x_{(n1+1)},\dots, x_{(n1+n2)}\rangle\circ f2, \dots,
\langle x_{n(k-1)+1},\dots, x_{n(k-1)+nk}\rangle\circ fk
>\circ f\tsk$$
The arguments of a derivative will be understood as {\it increments} $\trik x$
for argument $x$ of taken function.

The function's $n$-order derivative is defined as symmetric $n$-linear form. 
Therefore
it coincides with a symmetric polynome
\def\fract#1:#2 {{#1\over#2}}
$$\langle\trik_1, \trik_2,\dots,\trik_n\rangle\circ f^{(n)}=\fract 1:n!
\sum_{\sigma\in n!} \langle\trik_{\sigma(1)}, \trik_{\sigma(2)},\dots,
\trik_{\sigma(n)}\rangle\circ f^{(n)}\tsk$$
Each term of this  sum will be called a {\it concrete derivative graph}.
 The {\it virtual graph} is understood as a class of similar concrete 
derivative graphs. The
{\it symmetry number} $S$ of a virtual graph 
 helps to calculate the number of
similar members in such class. This number we shall call a 
{\it weihgt of virtual graph}.

We identify the  function's $n$-th order derivative with the  collection of all 
{\it virtual graphs} with
n  entrances. 
The virtual graph denotes the whole class of similar concrete derivative
graphs.

At the beginning we picture virtual graphs of the simple mapping 
for the first three derivatives
\vskip 12pt
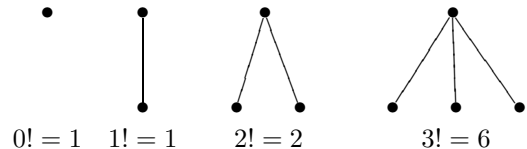
$$\matrix{
 \quad\xy ,+<\grynas,\aukstis>*{\bullet}**@{}\endxy \quad & 
\quad \xy *{\bullet},+<\grynas,\aukstis>*{\bullet}**\dir{-} \endxy\quad  &
\quad \xy @={*{\bullet},+<\zingsnis,\grynas>*{\bullet}} 
**@{}?<>(.5),
+<0cm,\aukstis>="virsus"*{\bullet}
@@{;"virsus"**@{}?(.95)**\dir{-}?>}\endxy \quad &
 \quad \xy @={*{\bullet},+<\zingsnis,\grynas>*{\bullet},+<\zingsnis,\grynas> 
*{\bullet}}
 **@{}?<>(.5),
+<0cm,\aukstis>="virsus"*{\bullet}
@@{;"virsus"**@{}?(.95)**\dir{-}?>}\endxy\quad \cr
0!=1 & 1!=1 &2!=2 &3!=6 \cr}$$

\vskip 12pt
Usually we don't picture the virtual graphs, it is enough that we can pick
a needed virtual graph and calculate its  weight number. The drawing of a
virtual graph and the calculation of its weight number can
be easily done automatically using a computer.

For the derivative of the simple mapping all members are similar, and we get
the virtual  graph identic with  concrete derivative graph, 
only  without ordering of
 argument increments
$$1=\fract n!:S \tsk$$
We shall choose a  concrete derivative graph with increasing
order $$\langle\trik_1, \trik_2, \dots, \trik_n \rangle$$
as {\it representing} for the virtual graph.

The derivative of some joint  mapping is found by changing the wedges of the
maping tree by the virtual graphs of derivatives for each  composed
simple function. Such changing must be correct: each increment of any
vertex must produce increment of elder vertexes, and is produced by
increment of some argument. At this nmoment such formulation will be
 sufficient for the drawing of virtual graphs. 

The number of similar members for each
virtual graph is calculated as a fraction  over graph's symmetry number
$$\fract n!:S  \tsk$$  
Each vertex $\ell$ has a symmetry number $S_\ell$ coinciding with the factorial
of the degree of the taken vertex
$$S_\ell=|\ell|!\tsk$$
The symmetry number of the whole graph is calculated as a product of
its vertices symmetry numbers
$$S=\prod S_\ell\tsk$$
We shall say that a virtual graph is {\it totally asymmetric} if the only
identical replacement $Id\in n!$ doesn't change the representing concrete
derivative graph.
Each replacement $\sigma\in n!$ will provide a different similar member, 
therefore 
the weight number for totally asymmetric virtual graph will be $n!$. 
It is hard  imagine such
possibility as interesting,
but it will be usual for more complex cases of derivative
calculation.

We shall say that a virtual graph is totally symmetric, if  every replacement
$\sigma\in n!$  provides the same representing concrete derivative graph.
In such case all similar
members coincide, and we shall have the unitary weight $1$ for totally symmetric
virtual graph.

We remark that for the linear function of two arguments $F(y,z)=y+z$,
the derivatives
for composition $F(f(x),g(x))$ will be calculated with usual binomial
coeficients
$$D^n=\sum_{0\mzl k\mzl n} \fract n!:k!(n-k)! f^{(k)} g^{(n-k)}\tsk$$
The virtual graph for such  function must be drawn whith 
coloured entrances and its weight
coincides with binomial coeficient. We can draw the virtual graph
 for the derivative of order $n=5$
with $k=2$ black entrances coresponding the increment of argument $f$
and $n-k=3$ white entrances corresponding the increment of argument $g$:

\vskip 12pt
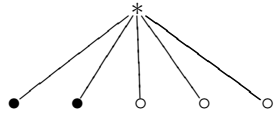
$$\matrix{
 \xy @={*{\bullet},+<\zingsnis,\grynas>*{\bullet},+<\zingsnis,\grynas> 
*{\circ},+<\zingsnis,\grynas>*{\circ},+<\zingsnis,\grynas>*{\circ},}
 **@{}?<>(.5),
+<0cm,\aukstis>="virsus"*{\ast}
@@{;"virsus"**@{}?(.93)**\dir{-}?>}\endxy\cr\noalign{\vskip 12 pt}
 S=2!\cdot 3!  \cr}$$
\vskip 12pt

{\bf Proposition 1.} The $n$-th order derivative 
of a joint function can be found
using the collection  of $n$-degree virtual graphs.
Their weights are equal to the cardinality of the whole symmetry group 
$n!$ divided by the
symmetry number $S$ of each  virtual graph.

{\bf Proof:} We must check that our procedure of virtual graph drawing
provides all needed concrete derivative graphs in the  derivative calculation, 
and only such
graphs. It is enough to check that we can get all concrete derivative
graphs, and then apply calculation of similar  concrete derivative 
graphs.

We apply the induction on the order $n$ of a joint function derivative.
For the zero order
derivative the proposition is trivial. We shall check the case $n=1$.
Calculating the first derivative,
all simple functions will be replaced by derivatives of first order.
One must be shure that we have got all  concrete derivative graph,
and every concrete derivative graph is provided by such procedure.

The induction step from $n$ to $n+1$: Taking the derivative of any
member presented by some concrete derivative graph, we differentiate some
vertex, and get additional entrance  of the concrete derivative graph of next 
degree
$(n+1)$. Also we must check that all concrete derivative graph of
degree $(n+1)$ can be obtained  in such manner from some concrete derivative
graph of degree $n$. It is done by  distraction anyone entrance
from the taken concrete derivative graph. 

\irodymopabaiga

{\bf 3. Ordering of virtual graphs.} For easy virtual graphs recognition we
need to choose simple ordering for all $n$-degree virtual graphs.
If this ordering will be useful for a wide class of users, it may
become standard.  We prefer to order all virtual graphs, and secondly
we induce this order to the set of virtual graphs of given degree $n$.
It will be called
a {\it natural order}. It is hard to imagine how somebody could
choose the best ordering only in
the set of $n$-degree graphs. 

We order the virtual derivative graphs lexicographically.
At first we order the
virtual derivative of the eldest function. If we have some of the eldest
functions, then we choose the order between them. In such sase we shall say 
that the top vertex is {\it ordered by the colour}. Each colour corresponds
to some sort of the eldest functions.

Then we order the vertices of one colour by the degree of this wertex.
We begin from the vertex of $0$-degree, and then go to the higher  degree.
If the degree is higher than $1$, then we order at first the vertex
from the left argument, and after we go to the next argument to the right.

The ordering of new vertices is the same: the colour, degree, the
younger vertex first from the left, and so on.
Therefore virtual graph will be represented by concrete graphs having
at left the younger graphs with
first colour and smaller degree.
\vskip 12pt
{\bf 4. The problems for the future.}
In geometrical calculations local coordinate change provides coordinate
change for various geometric objects. Such new change is calculated as
derivatives of joint function. The possible equality of joint functions
compels us to construct some virtual graphs.
In such a way we obtain new weights. 
They are obtained from earlier  weights with some concrete addition
operator. The earlier weights can  be called as {binomial} and they 
 present a free object for the
derivative calculation task. The question remains open, for what derivative
calculations such free object exists.

We shall give two examples of another free derivative calculation.
The first one provides the derivatives of inverse function $g=f\inv$.
Let these
functions operate over the points $x\in X$ and $y\in Y$
$$y=f(x)\kbl  x=g(y) \tsk$$
Then the derivatives of inverse function is calculated
$$Dg(y)=(Df(g(y))\inv \kbl$$
$$D^2g(y)= -\langle Dg(y),Dg(y)\rangle\dot D^2f(g(y))\dot Dg(y)\kbl $$
$$D^3g(y)=+3 \langle Dg(y),\langle Dg(y),Dg(y)\rangle\dot D^2f(g(y))\dot Dg(y)
\rangle \dot D^2f(g(y))\dot Dg(y)$$ 
$$-
<Dg(y),Dg(y),Dg(y)>\dot D^3f(g(y))\dot Dg(y)\tsk$$
The virtual derivative graphs are produced from the initial graph 
having only one
vertex and one wedge for the identic composition $1=f\circ g$.

\vskip 12pt
$$\xy*{\bullet} \ar@(ur,ul)c\endxy$$
\vskip 12pt

The virtual derivative graphs are drown without first order derivatives,
but these
derivatives must be written in final expression. The sign and weight of
virtual derivative graph are immediately appointed
in the same manner as in the previous case,
cl. Valiukevi\v cius [5], table V. We shall draw only virtual derivative
graphs having degree $n=3$. Under the graphs we shall write the graph's 
signed symmetry number. The first virtual graph will have a weight number 
$+3$ and second virtual graph will have a negative unitary numbe $-1$. A 
reader can verify that such weights are needed for the corresponding members
in the derivative formula.
\vskip 12pt
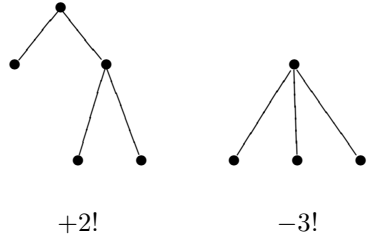
$$\matrix{%
\quad
\xy\save +<\grynas,\aukstis>**@{}*{\bullet}="virsune1"
\restore\POS +<\zingsnis,\grynas>;
@={;*{\bullet},+<\zingsnis,\grynas>*{\bullet}}
**@{}?<>(.5) +<\grynas,\aukstis> *{\bullet}="virsune2" 
@@{;"virsune2"**@{}?(.97)**\dir{-}?>}
@i \POS <0pt,0pt>;"virsune1":;@={"virsune1","virsune2"}
**@{}?<>(.5) +<\grynas,.6\aukstis> *{\bullet}="virsune" 
@@{;"virsune"**@{}?(.97)**\dir{-}?>}
\endxy\quad &
 \quad \xy @={*{\bullet},+<\zingsnis,\grynas>*{\bullet},+<\zingsnis,\grynas> 
*{\bullet}}
 **@{}?<>(.5),
+<0cm,\aukstis>="virsus"*{\bullet}
@@{;"virsus"**@{}?(.95)**\dir{-}?>}\endxy\quad \cr\noalign{\vskip 12 pt}
+2! & -3! \cr}$$
\vskip 12pt

The second example provides the derivatives of a solution of a
differential equation 
$$y'=f(y)\tsk$$
The following derivatives are obtained
$$y^{(2)}=f(y)\dot Df(y)\kbl$$
$$y^{(3)}=f(y)\dot Df(y)\dot Df(y)+ \langle f(y),f(y)\rangle\dot D^2f(y)\kbl$$
$$y^{(4)}=f(y)\dot Df(y)\dot Df(y)\dot Df(y)$$ 
$$+\langle f(y),f(y)\rangle\dot D^2f(y)\dot Df(y)$$ 
$$+ 3\langle f(y)\dot Df(y),f(y)\rangle\dot D^2f(y)$$
$$+\langle f(y),f(y),f(y)\rangle\dot D^3f(y)\tsk$$
The virtual derivative graphs are provided from the initial graph having only
one vertex representing the field $f$.

\vskip 12pt
$$\matrix{ \bullet \cr\noalign {\vskip 12pt}
f\cr}$$
\vskip 12pt

The virtual derivative graph is composed of
field derivatives $f^{(k)}$,
and each instance of derivative is provided by concrete
derivative graph with ordered vertices.

In this case the degree $n$ of virtual graph
is defined by cardinality of graph vertex set, and the "binomial"
weight is calculated with the earlier graph's symmetry number and additionally
with the new {\it complexity number}. We define the graph's cardinality 
as the number of its vertices. 
If some vertex $k$ has the
younger graphs with cardinality $c_1$, $c_2$, \dots, $c_{|k|}$ then the 
complexity number of taken vertex is defined by the product of these
cardinalities 
$$\tau_k=c_1\cdot\dots\cdot c_{|k|}\tsk$$
For the whole graph of $n$-degree 
the complexity number is equal to the product of
its vertices complexity numbers
$$\tau=\tau_1\cdot \dots\cdot \tau_n\tsk$$
For the $n$-degree virtual derivative graph the weight is calculated 
as $(n-1)!$
divided by the symmetry number and complexity number
$$ \fract (n-1)!:S{\tau} \tsk$$
The graphs are drawn in Valiukevi\v cius [5], table VII.
Now we shall draw only the virtual derivative graphs of degree $n=4$.
Under the graphs we shall write the graph's symmetry number and 
complexity number. All virtual graphs have the unitary weight numbers, 
only the third graph has the weight number $3$. A reader can verify that such 
weights are good for the corresponding members in the derivative formula.
$$\matrix{%
 \quad
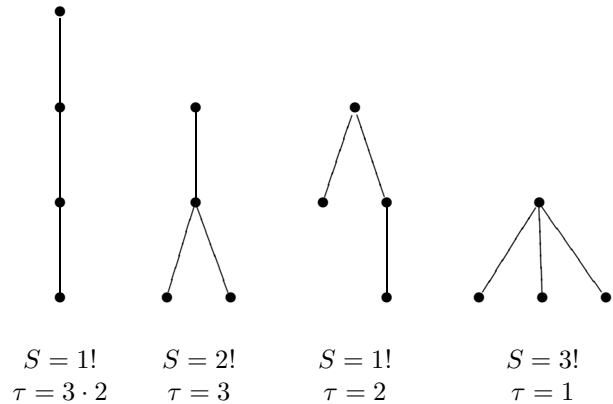
\xy *{\bullet},+<\grynas,\aukstis>*{\bullet}**@{-},
+<\grynas,\aukstis>*{\bullet}**@{-},+<\grynas,\aukstis>*{\bullet}**@{-} 
\endxy \quad 
& 
\quad \xy @={*{\bullet},+<\zingsnis,\grynas>*{\bullet}} 
**@{}?<>(.5)+<\grynas,\aukstis>="virsus";"virsus":*{\bullet}
+<\grynas,\aukstis>;**@{-};*{\bullet}
@@{;"virsus"**@{}?(.95)**\dir{-}?>}
\endxy \quad 
&
\quad \xy\save +<\grynas,\aukstis>
*{\bullet}="virsune1"
\restore\POS +<\zingsnis,\grynas>;
@={;*{\bullet}} 
**@{}?<>(.5) +<\grynas,\aukstis> *{\bullet}="virsune2" 
@@{;"virsune2"**@{}?(.97)**\dir{-}?>}
@i \POS <0pt,0pt>;"virsune1":;@={"virsune1","virsune2"}
**@{}?<>(.5) +<\grynas,\aukstis> *{\bullet}="virsune" 
@@{;"virsune"**@{}?(.92)**\dir{-}?>}
\endxy\quad
&
\quad \xy @={*{\bullet},+<\zingsnis,\grynas>*{\bullet},
+<\zingsnis,\grynas>*{\bullet}} 
**@{}?<>(.5),
+<0cm,\aukstis>="virsus"*{\bullet}
@@{;"virsus"**@{}?(.95)**\dir{-}?>}\endxy \quad
\cr
\noalign{\vskip 12pt}
S=1! &        S=2!   & S=1!   & S=3! \cr
\tau=3\cdot 2&\tau=3& \tau=2 & \tau=1 \cr}$$
\vskip 12pt

{\bf References:}

1. D. Gromoll, W. Klingenberg, W. Meyer 1968, 
Riemannische Geometrie im Gro\ss en,
Springer: Berlin / Heidelberg/ N. Y., 1971  Mir: Moskva.

2. A. Vinogradov, I. Krasil\v s\v cik, V. Ly\v cagin 1986, 
 Vvedenije v geometriju
nelinejnych differencialjnych uravnenij, Nauka: Moskva.

3. V. Arnold 1966,  Sur la g\'eom\'etrie diff\'erentielle des groupes de
Lie de dimension infinie et ses applications \`a l'hydrodynamique des fluides
parfaits, 319 - 361, Annales de l'Institut Fourier {\bf 16}: 1.

4. S. Lang 1962, Introduction to differentiable manifolds,
John Wiley \& Sons: London, 1967 Mir: Moskva.

5. G. Valiukevi\v cius 1979,
 Kombinatoriniai keliai i\v svestin\. ems skai\v ciuoti, 43 psl.+ I - XVIa 
(manuscript in Vilnius University Bibliotheca).

6. G. Valiukevi\v cius 1989,  On one measure contraction in lattice case, 
4 - 21,
Conspects on analysis, Vilnius.
\vskip 1cm

2011.01.21
\end{document}